\documentstyle[12pt,leqno]{article}
\def\af#1{\mbox{\bf\raise.086cm\hbox{\bf\tiny\bf/}\hskip-.14cm A$\!^{#1}$}}
\def\ar{\mbox{${\,\rightarrow\,}$}}
\def\ba#1{\begin{array}{#1}}\def\ea{\end{array}}
\def\be{\begin{equation}}\def\ee{\end{equation}}
\def\bi{\begin{itemize}}\def\ei{\end{itemize}}
\newtheorem{exx}{}[section]
\def\bex{\begin{exx}{\bf{Example. }}\em}\def\ex{\rm\end{exx}}

\def\bp{\vskip15pt\begin{exx}{\bf{Proposition. }}}\def\ep{\rm\end{exx}}
\def\bl{\vskip15pt\begin{exx}{\bf{Lemma. }}}\def\el{\rm\end{exx}}

\def\cc {\mbox{\hbox{\rm I\hskip-2pt K}}}
\def\dim{\hbox{dim\hskip2pt}}
\def\f{\raise2pt\hbox{$\varphi$}}
\def\ff#1{{\bf F}_{#1}}
\def\g{{\bf\Gamma}}
\def\gf{{\bf\Gamma}_{{\underline{f}}}}
\def\gm
{\mbox{$\raise1pt\hbox{\mbox{\footnotesize$\rm I$}}\!\!\!{\bf G}_m^
{\times n}$}}\def\gr#1{{\bf G}_{#1}}
\def\h{\hbox{\bf Hilb\hskip1pt}}
\def\id#1{\mbox{$\langle {#1} \rangle$}}
\def\ki{\mbox{\raise3pt\hbox{$\chi$}}}

\def\kp{\mbox{\large$\kappa$}}

\def\lar{\mbox{${\,\longrightarrow\,}$}}
\def\L{\mbox{${\cal L}$}}
\def\n{\noindent}
\def\O{\mbox{${\cal O}$}}
\def\pf{\vskip10pt\noindent{\bf Proof.\hskip10pt}}
\def\p#1{{\rm I}\!{\rm P}^{#1}}
\def\pd#1{\check{\rm I\! P}\hskip-8pt\phantom{\rm P}^{#1}}
\def\ps#1{\mbox{{\rm I$\!$P}{$\!$\rm\big(${#1}$\big)}}}
\def\qed{\phantom{|}\hfill\mbox{$\Box$}}
\def\rf#1{\hbox{\rm(\ref{#1})}}
\def\s{{\bf S}}
\def\so{\s^{nd}}
\def\st{\mbox{$\,|\,$}}
\def\t{{\bf T}}
\def\tkp{\wt{\kp}}
\def\wt{\widetilde}
\def\w{{\bf W}}
\def\z{{\bf Z}}
\def\wed#1{\mbox{$\stackrel{#1}{\wedge}$}}
\def\vez{\mbox{$\times$}}
\def\xi{\mbox{$\raise3pt\hbox{$\chi$}_
	{\raise-1pt\hbox{$_{\underline{f}}$
			}
	}
	     $}
	}

\begin{document}

\baselineskip=16pt
\parskip=7pt

\begin{center}{\bf Flatness of families induced\\
	by hypersurfaces on flag varieties}\vskip.5cm
		Israel Vainsencher
\footnote{$~ ^ {\heartsuit} ~$~ Partially supported by Brasil's CNPq.
Thanks are due to the UFMG for the warm atmosphere and hospitality during
the preparation of this manuscript.}$~ ^ {\heartsuit} ~$~\\
Departamento de Matem\'atica -- UFPE
\end{center}
\begin{abstract}
We answer a question posed by S. Kleiman concerning flatness of the
family of complete quadrics. We also show that any flat family of
hypersurfaces on Grassmann varieties induces a flat family of
intersections with the corresponding flag variety.
\end{abstract}
\section*{Introduction}

Let ~$\s$~ be the variety of complete quadrics, ~${\so}$~ the open
subset of nondegenerate quadrics and  ~${\ff{}}$~ the scheme of
complete flags in ~$\p n$. Let ~$\f_0:{\so}\ar \h({\ff{}})$~ be the
morphism that assigns to each  non\-degen\-erate quadric the locus of its
tangent flags. We prove the following.

{\bf Theorem. } \it ~$\f_0$~ extends to a morphism
{}~$\f:\s  \ar \h(\ff{})$\rm.

This answers affirmatively a question S. Kleiman asked in ([K], p.362).

We  first show that ~${\s}$~ parametrizes a flat family that
restricts, over  ~${\so}$, to the family of the graphs of the Gauss
map  (point $\mapsto$ tangent hyperplane) of nondegenerate quadric
hypersurfaces. The family pertinent to Kleiman's question is obtained
by composing the family of graphs with the appropriate flag bundle
(point $\in$ line $\subset \dots \subset $ hyperplane).

Our proof of flatness for the completed family of graphs relies on
Laksov's description [L] of Semple--Tyrrell's ``standard'' affine open
cover of \s.

The space of complete conics has recently reappeared as a simple
instance of Kontsevich's spaces of stable maps (cf. [P]). It is also
instrumental for the counting of rational curves on a K3 surface double
cover of the plane (cf. [V]). Complete quadric surfaces play a role
in Narasimhan--Trautmann [NT] study of a compactification of a space of
instanton bundles.

We also show that \em any flat family of hypersurfaces on Grassmann
varieties induces a flat family of subschemes of the corresponding
flag variety \rm (cf. \rf{prop}).

This statement was first obtained as an earlier attempt to answer
Kleiman's question. We observe that for the case of quadric
hypersurfaces the family described in the proposition  does not induce
the family  of tangent flags. In fact, for conics it yields a double
structure on the graph of the Gauss map.(cf. \S\ref{fim} for details).

\section{The tangent flag to a smooth quadric}

 Write $x =(x_1,\dots,x_{n+1})$ (resp. $y =(y_1,\dots,y_{n+1})$) for
the vector of homogeneous coordinates in $\p n$ (resp. $\pd n$ ). Let
{}~$\ff{0,n-1}\subset\p n\times\pd n$~ be the incidence correspondence
``point $\in$ hyperplane''. It is the zeros of the incidence section
{}~$x\cdot{}y$~ of $\O_{\p n}(1)\otimes\O_{\pd{\,\, n}}(1)$.

Let ~$\kp\subset\p n$~ denote a smooth quadric represented by a
symmetric matrix $a$. The Gauss map ~$\gamma:\kp\ar\pd n$~  is given
by $x\mapsto y=x\cdot a$. Hence we have
$$
\gamma^*(\O_{\pd{\,\, n}}(1))=\O_{\p n}(1)_{|\kp}.
$$
The tangent flag ~$\tkp\subset\ff n$~ of ~$\kp$~ is equal to the
restriction of the flag bundle
$$
\ff n\ar\ff{0,n-1}\subset\p n\times\pd n \vspace{-1pt}
$$
over the graph ~$\g_{\mbox{$\kappa$}}$~ of ~$\gamma$. Consequently,
flatness of the family ~$\{\tkp\}$~ of tangent flags
is equivalent to flatness of the
family of graphs ~$\{\g_{\mbox{$\kappa$}}\}$. The latter will be
handled in \S\ref{graphs}.

We proceed to compute the Hilbert polynomial of the graph $\g$ of the
Gauss map of a general quadric hypersurface ~$\kp\subset\p n$.
\bl
Notation as above, the Hilbert polynomial
{}~$\ki\big(\O_{\g}(\L^{\otimes t})\big)$~ with respect to
$$
\L=\big(\O_{\p n}(1)\otimes{}\O_{\pd{\, n}}(1)\big)_{|\g}
$$
is equal to
$$
{2\,t+n\choose n}-{2(t-1)+n\choose n}.
$$
\label{hilbgen}\el

\vspace{-15pt}\pf
We have $\L\cong\O_{\p n}(2)_{|\kp}$ under the identification $\g\cong\kp$.
Thus we may compute
$$
\ba{cl}\ki(\L^{\otimes t})
&= \ki\big(\O_{\p n}(2t)\big)_{|\kp}\\
\noalign{\vskip5pt}
&=\ki\big(\O_{\p n}(2t)\big) - \ki\big(\O_{\p n}(2t-2)\big)\\
\noalign{\vskip7pt}
&={2\,t+n\choose n}-{2(t-1)+n\choose n}.
\ea
$$\vspace{-30pt}\vskip-7pt\qed

\section{Hilbert polynomial of loci of rank 1 matrices}

The image of the Segre imbedding $\p n\times\p n\ar\p N$ is the
variety of matrices of rank one. The image $\Delta$ of the diagonal
$\p n \ar \p n \times\p n\ar\p N$ is the subvariety of $symmetric$
matrices of rank one. It's Hilbert polynomial is easily found to be
$$
\dim(H^0(\Delta,\O_{\p N}(t))~=~{2t+n\choose n}.
$$
The bi-homogeneous ideal $I_{\Delta}$ of the diagonal is generated by
the 2$\times2$ minors of the matrix
\be\label{2x2}
\left [\ba{cccc}
x_{{1}}&x_{{2}}&\dots&x_{{n+1}}\\\noalign{\medskip}
y_{{1}}&y_{{2}}&\dots&y_{{n+1}}
\ea\right].
\ee
Write
$$
S=k[x_1,\dots, x_{{n+1}},y_1,\dots ,y_{{n+1}}]
$$
for the polynomial ring in $2n+2$ variables, and let $S_{i,j}$ denote
the space of bi-hom\-oge\-neous polynomials of bi-degree $(i,j)$. We
have
$$
\dim_k\,S_{t,t}\big/(I_{\Delta})_{t,t}~=~{2t+n\choose n}.
$$
Quite generally, for a closed subscheme $X\subseteq \p m\times\p n$
defined by a bi-homogeneous ideal $I\subseteq S$ we have (cf. [KTB],
p. 189)
$$
H^0\big(X,\O_{\p m}(t)\otimes{}\O_{\p { n}}(t)_{|X}\big)~=~
S_{t,t}\big/(I)_{t,t}~~\hbox{ for all }~t >>0.
$$
Indeed, the homomorphism
\be\label{rs}\ba{ccc}
R=k[z_{i,j}]\big/\id{z_{i,j}z_{r,s}-z_{r,j}z_{i,s}}&\lar &
		S\\
z_{i,j}&\mapstochar\longrightarrow& x_iy_j
\ea\ee
maps $R_t$ isomorphically onto $S_{t,t}$. Let the bi-homogeneous ideal
$I\subseteq S$ be generated by polynomials of some fixed bidegree
($d,d$). Its inverse image via \rf{rs} generates a homogeneous ideal
$I^\#\subseteq R$.
We have
$$
(R/I^\#)_t\wt{\lar}S_{t,t}/I_{t,t}.
$$
Now for $t>>0$ we may write
$$
(R/I^\#)_t~=~H^0\big(X,\O_{\p N}(t)_{|X}\big)~=~
H^0\big(X,\O_{\p m}(t)\otimes{}\O_{\p {n}}(t)_{|X}\big).
$$

Let $L(I)$ denote the monomial ideal of initial terms of $I$ with
respect to some bi-graded monomial order. Then we have the equality of
Hilbert functions,
$$
\f_I(i,j)~=~\f_{L(I)}(i,j).
$$
This is rather standard: let $f_1,\dots,f_k$ be linearly independent
forms of bidegree $(i,j)$ in $I$. Replacing if needed each $f_\mu$ by
$f_\mu-cf_\nu$ for suitable $c\in k$, we may assume their initial
terms $L\,f_\mu \ne L\,f_\nu $. Hence the initial terms $Lf_1, \dots,
Lf_k$ are linearly independent monomials in $L(I)_{i,j}$. This shows
that $ \f_I(i,j) \leq \f_{L(I)}(i,j). $ Conversely, pick monomials
\,$g_1> \dots> g_k$\, in $L(I)_{i,j}$. We have each $g_\mu = Lf_\mu$
for some  $f_\mu \in I_{i,j}$. It follows that $f_1, \dots, f_k$ are
linearly independent.

\bl\label{k0}
Let ~$\g_0$~ be the subscheme of $\p n\times\pd n$ defined by the ideal
$$
\id{x_iy_j\st 1\leq i<j\leq n+1}+\id{\sum x_iy_i}.
$$Then we have
$$
\f_{\g_0}(t)={2\,t+n\choose n}-{2(t-1)+n\choose n}.
$$
\el
\pf
The whole point is to notice that the $ x_iy_j$ span the ideal of
initial terms of $I_\Delta$ with respect to a suitable order.\footnote
{I'm indebted to P. Gimenez for his precious help on this
matter.} ~In fact, the set of $2\times2$ minors of \rf{2x2} is known
to be a (universal) Gr\"obner basis for $I_\Delta$ (see Sturmfel
[BS]). By the above discussion, we may write
$$
\f_{L(I_\Delta)}(t)=\f_{I_\Delta}(t)={2t+n\choose n}.
$$
One checks at once that $\sum x_iy_i$ is a nonzero divisor mod
$L(I_\Delta)$ (see \rf{rmk} (i)). Therefore
$$
\f_{\g_0}(t)=\f_{L(I_\Delta)}(t)-\f_{L(I_\Delta)}(t-1).\vspace{-10pt}
$$\qed

We will deduce flatness for the ``completed'' family of Gauss maps
from the fact that the above Hilbert polynomial at the special point
{}~$\g_0$~ coincides with the generic one \rf{hilbgen}.

\section{Semple-Tyrrell-Laksov cover}

Let ${\bf U}_n$ denote the group of lower triangular unipotent
($n+1$)-matrices. Thus, ${\bf U}_n$ is isomorphic to the affine space
$\af{n(n+1)/2}$ with coordinate functions ~$u_{i,j},\,1\leq j\leq
i-1,\,i=2\dots n+1$. These are thought of as entries of the matrix,
$$
u~=~\left [\begin {array}{lllll} 1&0&0&\cdots&0\\
\noalign{\medskip}u_{{2,1}}&1&0&\cdots&0\\
\noalign{\medskip}u_{{3,1}}&u_{{3,2}}&1&\cdots&0\\
\noalign{\medskip}\cdots&\cdots&\cdots&\cdots&\cdots\\
u_{{n+1,1}}&u_{{n+1,2}}&u_{{n+1,3}}&u_{{n+1,n}}&1
\end {array}\right ].
$$
Let $d_1,\dots,d_n$ be coordinate functions in $\af n$.
Put
\be\label{d1}
d^{(1)}=\left [\begin {array}{lllll} 1&0&0&\cdots &0\\
\noalign{\medskip}0&d_{{1}}&0&\cdots &0\\
\noalign{\medskip}0&0&d_{{1}}d_{{2}}&\cdots &0\\
\noalign{\medskip}\cdots&\cdots&\cdots&\cdots&\cdots\\
0&0&0&\cdots &d_{{1}}d_{{2}}\cdots d_{{n}}
\end {array}\right ].
\ee
For a matrix $A$ let it's $i$th adjugate be the matrix $\wed i A$ of
all $i\times i$ minors. We denote by $d^{(i)} $ the matrix
obtained from $\wed i d^{(1)}$ by removing the common factor
$d_1^{i-1}d_2^{i-2}\cdots d_{i-1}$. $E.g.,$ for $n=3$ we have
$$\ba{ll}
d^{(1)} &= {\rm diag}(1,\,d_1,\,d_1d_2,\, d_1d_2 d_3)\\
\noalign{\medskip}d^{(2)} &= {\rm
diag}(d_1,\,d_1d_2,\,d_1d_2d_3,\,d_1^2d_2,\,d_1^2d_2d_3,\,d_1^2d_2^2d_3)
\big/(d_1)
\\
&= {\rm diag}(1,\, d_2,\, d_2d_3,\, d_1d_2,\, d_1d_2d_3,\, d_1d_2^2d_3)
\\\noalign{\medskip}
d^{(3)} &=
{\rm diag}(1                    ,\, d_3     ,\,d_2d_3   ,\,d_1d_2d_3).
\ea$$
The map ~$ {\bf U}_n\times\af n\ar \s\subset{
\prod_{i=1}^{i=n}\ps{
S_2(\bigwedge\hskip-8pt\raise6pt\hbox{$^i$}\,\,k^{n+1*})}}$~
defined by sending $(u,\,d)$ to
$$
\big(u\,d^{(1)}\,u^t,\,(\wed2\,u)\,d^{(2)}\,\wed2u^t,\dots,(\wed nu)\,d^{(n)}\,
\wed nu^t\big)
$$
is an isomorphism onto an affine open subset $\s^0$ of \s.
The variety of complete quadrics may be covered by translates of
$\s^0$ (cf. Laksov [L]).

Let  $\s^0_d \cong{\bf U}_n\times\af n\!_{d}$ be the
principal open piece defined by $d_1d_2\cdots d_n\ne0$. It maps
isomorphically onto an open subvariety  of  $\so$.

\section{Graph of the Gauss map}\label{graphs}

The variety ~$\so$~ of nondegenerate quadrics parametrizes a flat
family of graphs of Gauss maps. For a nondegenerate quadric
represented by a symmetric matrix $a\in\so$ the Gauss map is given by
$x\mapsto y=x\cdot a$.  We define $\cc^{nd}\subset\so\vez\p n\vez\pd
n$ by the bi-homogeneous ideal generated by the incidence relation
{}~$x\cdot y$~ together with the 2$\times$2 minors of the 2$\times(n+1)$
matrix with rows ~$y,\,x\cdot z$,~ where \,$z$\, denotes the generic
symmetric matrix. Clearly ~$\cc^{nd}\ar\so$ is a map of ${\bf
GL}_{n+1}-$homogeneous spaces.

Now write $a = vc^{(1)}v^t$ with $v\in{\bf U}_n,\, c\in\af n\!_d$
\,($c^{(1)}$\, as in \rf{d1}), and put $x'=xv$, $y'=y(v^{-1})^t$. We
have $y=xa$ iff $y'= x'c^{(1)}$. Let
\be
{\cc}^0_d ~\subset~ \s^0_d \times\p n\times\pd n.
\label{g0d}
\ee

be defined by  ~$x\cdot y$~ together with the 2$\times$2 minors of the
2$\times(n+1)$ matrix
\be\label{x'd}
\left [\begin {array}{ccccc} x'_{{1}}&d_1x'_{{2}}&d_1d_2x'_{{3}}&\dots&
d_1\cdots d_nx'_{{n+1}}\\
\noalign{\medskip}y'_{{1}}&y'_{{2}}&y'_{{3}}&\dots&y'_{{n+1}}
\end {array}\right ]
\ee
where we put $x'_j=\sum_iu_{ij}x_i$ and likewise $y'_j$ denotes the
$j$th entry of $y(u^{-1})^t$. Thus ${\cc}^0_d$ is the
total space of the family of Gauss maps parametrized by $\s^0_d$.
Note that ${\cc}^0_d \ar \s^0_d$ is a smooth quadric
bundle. Its fibre over $(I,(1,\dots,1)) \in {\bf U}_n \times \af
n_{d}$ is equal to the quadric given by \,$\sum x_i^2$\, inside the
``diagonal'' \,$y_1=x_1, \dots,y_{{n+1}}= x_{{n+1}}$\, of ~$\p
n\times\pd n.$

Let
\be\label{gg0} {\cc}^0  ~\subset~ \s^0 \times\p n\times\pd n
\ee
be defined by  ~$x\cdot y$~ together with the ideal
\be\label{J}
\ba{cl}J=&\langle{}x'_1y'_2-d_1y'_1x'_2,\dots,\,x'_1y'_{n+1}-
d_1\cdots{}d_ny'_1x'_{n+1},\\
	 &\,\,\,x'_2y'_3-d_2y'_2x'_3,\dots,\,x'_ny'_{n+1}-d_ny'_nx'_{n+1}
\rangle
\ea\ee
obtained by cancelling all $d_i$ factors occurring in the above
2$\times$2 minors. We obviously have
{}~$ {\cc}^0_{\,|\,\s^0_d} = {\cc}^0_d$.

We will show that ~$ {\cc}^0$~ is the scheme theoretic closure of ~$
{\cc}^0_d$ in $\s^0 \vez \p n\vez\pd n$ (cf. \rf{clos}).

\section{A torus action}

Notation as in \rf{d1}, imbed  ~\gm~ in ${\bf GL}_{n+1}$ by sending
$c=(c_1,\dots,c_n)\in\gm$ to $c^{(1)}={\rm diag}(1,\,c_1,\,c_1c_2,\dots)$.
We let  ~\gm~ act on $\s^0$ by
$$
c\cdot (v,b) = (c^{(1)}\, v\, (c^{(1)})^{-1},~(c_1^2b_1, \dots c_n^2b_n)).
$$
This action is compatible with the natural action of ${\bf GL}_{n+1}$
on the space $\ps{S_2(k^{n+1*})}$ of quadrics, $i.e.,$ for a symmetric
matrix ~$a(v,b)\, :=\, v\, b^{(1)}\, v^t$~ as above, we have
$$
\ba{cl}
c^{(1)}\cdot{}a(v,b)&=c^{(1)}\, a(v,b)\,(c^{(1)})^t =
		c^{(1)} \,v\,b^{(1)}\,v^t \,(c^{(1)})^t
\\\noalign{\medskip}&=c^{(1)}\,v\,(c^{(1)})^{-1}\,c^{(1)}\,b^{(1)}
		\,c^{(1)}\,((c^{(1)})^t)^{-1}\,v^t\, (c^{(1)})^t
\\\noalign{\medskip}&=c^{(1)}\,v\,(c^{(1)})^{-1}\,(c^{(1)})^2\,b^{(1)}
		\,((c^{(1)})^t)^{-1}\,v^t\,(c^{(1)})^t
\\\noalign{\medskip}&=a(c\cdot(v,b))\,.
\ea$$
It can be also easily checked that ~\gm~ acts compatibly on
$\s^0\times\p n\times\pd n$ and ${\cc}^0$ is invariant.
Indeed, let $((v,b),x,y)\in{\cc}^0$. Pick $c\in\gm$.
We have
$$
c\cdot((v,b),x,y)=( (c^{(1)}\, v\, (c^{(1)})^{-1},~(c_1^2b_1, \dots
c_n^2b_n)),~x\,(c^{(1)})^{-1}, \,y\,(c^{(1)})^{t}).
$$
Now $x'=xv$ changes to
$$
x'' ~=~ (x\, (c^{(1)})^{-1})\, (c^{(1)}\, v\, (c^{(1)})^{-1}) ~=~
x\,v\, (c^{(1)})^{-1} ~=~ x'\, (c^{(1)})^{-1}
$$
so that the first row $x'\,b^{(1)}$ in \rf{x'd}
(evaluated at $((v,b),x,y)$)
changes to
$$
x''\,(b^{(1)}\,(c^{(1)})^2) ~=~ x'\, (c^{(1)})^{-1}
\,(b^{(1)}\,(c^{(1)})^2) ~=~  x'\, (b^{(1)}\, c^{(1)}).
$$
Similarly,
$y' ~=~ y\,(v^{-1})^t$ changes to
$$
y'' ~=~ (y\,(c^{(1)})^{t})\,
((c^{(1)}\, v\, (c^{(1)})^{-1})^{-1})^t  ~=~
y\,(v^{-1})^t\,(c^{(1)})^{t})
 ~=~ y'\, c^{(1)}.
$$
Therefore \rf{x'd} changes to the matrix
with rows $x'\, (b^{(1)}\, c^{(1)})$ and $y'\, c^{(1)}$.
Thus evaluation of \rf{J} at $c\cdot((v,b),x,y)$ and at $((v,b),x,y)$
differ only by nonzero multiples.

\bl
The orbit of $(I,0)\in\s^0$ is the unique closed orbit.
\el
\pf
Conjugation of \,$v\in {\bf U}_n$\, by the diagonal matrix \,$c^{(1)}$\,
replaces each entry \,$v_{ij},\,j<i$\, by
$$
\ba{cl}(c^{(1)}\,v\,(c^{(1)})^{-1})_{ij} &=
c^{(1)}_{ii}\,(v\,(c^{(1)})^{-1})_{ij}
 =  c^{(1)}_{ii}\, v_{ij}\,((c^{(1)})^{-1})_{jj}\\
	&= v_{ij}\,c^{(1)}_{ii}/c^{(1)}_{jj}= v_{ij}\,\,c_{i-1}\cdots{}c_{j}.
\ea$$
Thus, letting \,$c\ar0$, we see that ($I,0$) is in the orbit closure
{}~$\overline{\gm\cdot(v,b)}$

\vspace{-10pt}\qed

\section{Proof of the theorem}

\bl\label{cc0}
Notation as in \rf{gg0} , the family ~${\cc}^0\ar\s^0$~ is
flat.
\el
\pf
Since ~${\cc}^0 \ar \s^0$~ is equivariant for the ~$\gm-$action, it
suffices to check that the Hilbert polynomial of the fiber over the
representative $(I,0)$ of the unique closed orbit is right, $i.e.,$
coincides with the generic one (cf.  Hartshorne [H], thm.9.9, p.261).
Evaluating \rf{J} at $(I,0)$ yields the monomial ideal \rf{2x2}. We
are done by virtue of \rf{hilbgen} and \rf{k0}.\vspace{-25pt}

\qed

\bl
Let $f:X\ar Y$ be a flat, surjective morphism of schemes. If $U\subseteq Y$
is open and schematically dense in $Y$ then  $f^{-1}U$ is open and
schematically dense in $X.$
\el
\pf
We may assume $X,\,Y$ affine. Let $A\subseteq B$ be a flat ring
extension and let $a\in A$ be such that ${\rm Spec}\,A_a$ is
schematically dense in ${\rm Spec}\,A$. This means that every element
in ker\,($A\ar A_a$) is nilpotent. Flatness implies ker\,($B\ar
B_a$)=ker\,($A\ar A_a)\otimes{}B$. Hence ${\rm Spec}\,B_a$ is
schematically dense in ${\rm Spec}\,B$. \vspace{-25pt}

\qed

\bl\label{clos}
Notation as in \rf{gg0} and \rf{g0d}, we have that ~$\cc^0$~ is
equal to the scheme theoretic closure of ~$\cc^0_d$.
\el
\pf
In view of \rf{cc0}, we may
apply the previous lemma to ~$\cc^0\ar\s^0\supset\s^0_d.$
\vspace{-25pt}

\qed

\bl
Let $G$ be an algebraic group and let
$$
\ba{ccc}X^0 &\subset &X \\
\downarrow~&&\downarrow\\
Y^0&\subset &Y
\ea$$
be a commutative diagram of maps of ~$G-$varieties. Let
\,$\overline{X}$,~  $\overline{Y}$  denote the closures of
$X^0,\,Y^0.$  If  ~$\overline{X}\ar\overline{Y}$~ is flat over a
neighborhood of a point in each closed orbit then
{}~$\overline{X}\ar\overline{Y}$~ is flat.
\el
\pf Immediate.
\qed

We may now finish the proof of the theorem. Let ~$\cc\subset\s\vez\p
n\vez\pd n$~ be the scheme theoretic closure of ~$\cc^0$. We have
{}~$\cc\cap\big(\s^0\vez\p n\vez\pd n\big) = \cc^0$ flat over ~$\s^0$~
by \rf{cc0}. The latter is a neighborhood of a point in the unique closed
orbit of ~$\s$. Now apply the previous lemma to ~$G={\bf GL}_{n+1}$,
$X=\s\vez\p n\vez\pd n$, $Y=\s$, $Y^0=\s^{nd}$, $X^0=\cc^{nd}$.
Finally, since the family of tangent flags is defined by the fibre square,
$$
\ba{ccc}
\wt{\cc}&	\lar	&\ff n \, \vez \,\s \\
\downarrow&		&\downarrow\\
\cc	&	\lar	&\ff{0,n-1}\,\vez\,\s
\ea
$$
the composition ~$\wt{\cc}\ar\cc\ar\s$~ is flat.

\section{Final remarks}\label{fim}

\begin{exx}\em \label{rmk}
(i) The primary decomposition of the monomial ideal in \rf{k0} can be
checked to be given by
$$
 \id{x_1,\,x_2,\dots, x_n} \cap \cdots \cap
\id{x_1,\dots, x_i,\,y_{i+2}, \dots, y_{n+1}} \cap  \cdots
\cap \id{y_2,\,y_3, \dots, y_{n+1}}.
$$
Thus enlarging it to include the nonzero divisor \,$x\cdot{}y$\,
we see that the special fiber \,$\g_0$ \, presents no imbedded component.

(ii) The example of $\p n$ acted on by the stabilizer of a point,
blown up at that point might clarify why we were not able
to show directly that the closure of ~$\cc^{nd}$~ is flat over ~$\s$.

(iii) For $n=1$ we may write the following global equations
for $\cc$. Let $z,\,w$ be a pair of symmetric 3\vez3 matrices of
independent indeterminates. Then  ~$\cc \subset \p5 \vez \pd5 \vez \p2
\vez \pd2$~ is given by the 2\vez2 minors of the 2\vez3 matrices with
rows $x\cdot z,\,y$ and $x,\,y\cdot w$, in addition to the incidence
relation $x\cdot y$ together with the equation
{}~$3z\cdot w={\rm trace}(z\cdot{}w)I$~for $\s\subset\p5 \vez \pd5$.
It would be nice to give a similar description for higher dimension.

(iv) Still assuming $n=1$,
put
$$
\g=\{(P,\,\ell,\,\kp,\,\kp')\in\p2\vez\pd2\vez\p5\vez\pd5\st
P\in\kp\cap\ell,\,\ell\in\kp'\}.
$$
It is easy to check that ~$\g_{|\s}=\cc$~ as sets.
Furthermore, $\g$ may be endowed with a natural scheme structure such
that $\g\ar\p5\vez\pd5$ is flat and with Hilbert polynomial of its
fibers equal to $4t$. Thus, $\g_{|\s}\ar\s$ is a family of double
structures of genus one on the fibers of $\cc$.
\end{exx}

In fact, we have the following.
\bp\label{prop}
 Any flat family of hypersurfaces on Grassmann varieties induces a
flat family of subschemes of the corresponding flag variety.
\ep

Before considering the general case,
we describe the situation in the projective plane. Thus, let
$$
{\ff2}\subset\p2\times\pd2
$$
be the incidence correspondence ``point ~$\in$~ line''.
Let ~$f_0$~ (resp. ~$f_1$) denote a curve in ~$\p2$~ (resp. ~$\pd2$).
Set
$$
\gf:=(f_0\times f_1)\cap\ff2.
$$
\n Then ~$\gf$~ is easily seen to be regularly imbedded of codimension
2  in ~${\ff2}$~ (cf. \rf{claim}). Moreover, its Hilbert polynomial
 with respect to the ample sheaf ${\cal
O}_{\p2}(1) \otimes {\cal O}_{\pd{\,\,2}}(1)$ restricted to ~${\ff2}$~
depends only on the degrees, say ~$d_0,\,d_1$~ of ~$f_0,f_1$. In fact,
the Koszul complex that resolves the ideal of ~$f_0\times f_1$~ in
$\p2 \times \pd2$~ restricts to a resolution of ~$\gf$~ in ~${\ff2}$.
One finds the Hilbert polynomial,
\be
\label{xi}
\xi\hskip-2pt(t)=(d_0+d_1)t- d_0 d_1(d_0+d_1-4)/2 .
\ee
Therefore, as in the final argument for the proof of \rf{cc0},
the parameter space of pairs ~$(f_0,f_1)$, call it ~$\t$
(=$\p{n_0}\times\p{n_1}$~ for suitable ~$n_0,n_1)$,
carries  a flat family of curves on ~${\ff2}$.\, Precisely, let
$$
\w\!_0\subset \p2\times\p{n_0}\hbox{ and }
\w\!_1\subset \pd2\times\p{n_1}
$$
denote the total spaces of the universal plane curve parametrized by
$\p{n_i}$. Then
$$
\g:=(\w\!_0\times\hskip-.38cm\raise-.25cm\hbox{$_{\p2}$}
\w\!_1)\bigcap{\ff2}\longrightarrow \t
$$
is a flat family of curves in ~${\ff2}$, with fiber ~$\gf$.

\vskip10pt
For the proof of \rf{prop} we
let ~$\gr{r,n}$~ denote the grassmannian of projective subspaces
of dimension \,$r$\, of ~$\p n$.

Recall that the dimension of the variety of complete flags
$\ff n\subset\prod\gr{i,n}$~ is
$$
\dim{\ff n}=1+\cdots+n.
$$

The proposition is an easy consequence of the following.

\bl \label{claim}
 Let ~$f_0,f_1,\dots,f_n$~ be arbitrary hypersurfaces of
points, lines, \dots, hyperplanes in the appropriate grassmannians of
subspaces of ~$\p{n+1}$. Then the intersection
$$
\gf:= (f_0\times\cdots\times f_n)\cap\ff{n+1}
$$
is of codimension ~$n+1$~ in ~${\ff{n+1}}$.
\ep

\pf
We shall argue by induction on ~$n$.

We may assume all ~$f_i$~ irreducible. Let ~$n=1$. Pick a line ~$h\in
f_1$. Set
$$
h^ {(0)}=\{P\in\p2\st{}P\in{}h\}.
$$
The fiber ~$(\gf)_{h} \simeq
h^ {(0)}\cap f_0$~ is zero dimensional unless ~$h^{(0)}=f_0$. This
occurs for at most one ~$h\in{}f_1$, hence ~$\gf$~ is 1--dimensional
(otherwise most of its fibres over ~$f_1$~ would be at least
1--dimensional).

For the inductive step, we set for ~$h\in \pd{n+1}$,
\be\label{hr}
h^ {(r)}=\{g\in\gr{r,n+1}|g\subseteq h\}\simeq\gr{r,n}.
\ee
If the
intersection
$$
f'_r  =h^{(r)}\cap f_r
$$
were proper for all ~$r$
and ~$h\in f_n$~ then we would be done by induction. Indeed, we have
$$
(\gf)_h \simeq(f'_0 \times\cdots\times f'_{n-1})\cap\ff n.
$$
By the induction hypothesis, this is of the right dimension
$$
1+\cdots +n-n =1+\cdots +(n-1).
$$
Since ~$h$~ varies in the ~$n-$dimensional hypersurface ~$f_n$~ of
$\gr{n,n+1}=\pd{n+1}$, we would have
$$
\dim \gf=\big(1+\cdots +(n-1)\big)+n=\big(1+\cdots +(n+1)\big)
-(n+1)
$$
as desired.

However, just as in the case ~$n=1$,  it may well happen that the
intersection ~$h^{(r)}\cap f_r$~ be \it not \rm proper for some ~$h,r$.
Thus it remains to be shown that, whenever \dim$(\gf)_h$~ exceeds the
right dimension, say by ~$\delta$, the hyperplane ~$h$~ is restricted to
vary in a locus of codimension at least ~$\delta$~ in ~$f_n$. This is
taken care of by the lemma below.

\bl\label{lema}
Notation as in \rf{hr},
for ~$r=0,\dots,n$~ we have
$$
\dim \{ h\in\pd{n+1}\ |\ h^{(r)}\subseteq\ f_r\} \leq r.
$$
\el
\n {\it Proof}. Let ~$\ff{r,n}\subset\pd{n+1}\times\gr{r,n+1}$~ be the
partial flag variety. Form the diagram with natural projections,
$$
\ba{lcr}&\ff{r,n}&\\\noalign{\vskip5pt}
\raise6pt\hbox{$\pi_n$}\! \mbox{\huge$\swarrow$} & & \mbox{\huge$\searrow$}
\raise6pt\hbox{$\pi_r$}~~~~~\\
\noalign{\vskip3pt}\pd{n+1} && ~~\gr{r,n+1}
\ea$$
For ~$g_r\in\gr{r,n+1}$, set
$$
g_r^ {(n)}=\{h\in\pd{n+1}\ |\ g_r\subseteq h\}.
$$
We have ~$g_r^ {(n)} \simeq\p{n-r}$~ whence it hits any
subvariety of ~$\pd {n+1}$~  of dimension ~$\geq r+1$. In other
words, for any subvariety ~$\z\subseteq \pd{n+1}$~ such that
\dim ~$\z\geq r+1$, we have
$$
\ba{ccl}\pi_r\pi^ {-1}_n\z    &=&
	\{g_r\ |\ \exists\, h\in \z\hbox{ s.t. } h\supseteq g_r\}\\
		&=&     \{g_r\ |\ g_r^{(n)}\cap\z\ne\emptyset\}\\
	&=&\gr{r,n+1}.
\ea$$
The lemma follows by taking \z$=
\{ h\in\pd{n+1}\ |\ h^{(r)}\subseteq\ f_r\}$.
Indeed, if \dim\z$\geq r+1$, then for all ~$g_r\in\gr{r,n+1}$~ there exists
$h\in\z\hbox{ s.t. } h\supseteq g_r$, so ~$g_r\in h^{(r)}\subseteq f_r$,
contradicting that ~$f_r$~ is a hypersurface of ~$\gr{r,n+1}$. \qed(for
\rf{lema})

Continuing the proof of \rf{claim} we consider the stratification of
$f_n$~ by the condition of improper intersection of ~$f_r$~ with ~$h^
{(r)}$, namely,
$$
\matrix {f_{n,0}&=&\{h\in f_n\ |\ h^{(0)}\subseteq f_0\},&\cr
f_{n,1}&=&\{h\in f_n\ |\ h^{(1)}\subseteq f_1\}&\hskip-.25cm\setminus&
\hskip-.75cm f_{n,0},
\cr
&\vdots&&&\cr
f_{n,n}&=&\{h\in f_n\ |\ h^ {(n)}\subseteq f_n\}&\hskip-.25cm
\setminus&\hskip-.15cm \bigcup\hskip-.45cm\raise-.25cm\hbox{$_{j<n}$}f_{n,j}.}
$$
We will be done if we show
$$
\dim(\gf)_h\ \leq 1+\cdots+n-
r\quad\forall\ h \in f_{n,r}.
$$
We have already seen that ~$\dim(\gf)_h=1+\cdots +n-1$~ for ~$h$~ in
$f_{n,n}$. Also, for ~$r=0$, the desired estimate holds because we have
$(\gf)_h\subseteq (\ff{n+1})_h\simeq\ff n$~ and ~$\dim \ff
n=1+\cdots+n.$~  Let ~$r>0$~ and pick a hyperplane ~$h \in f_{n,r}$. Then
the intersections,
$$
f'_i = h^ {(i)}\cap f_i,
$$
are proper for ~$i=0,\dots,r-1$, whereas for the subsequent index, we
have
$$
h^ {(r)}\cap f_r = h^ {(r)}\simeq \gr{r,n}.\vspace{-10pt}
$$
Thus, we may write,
$$
(\gf)_h\hookrightarrow
\big(f'_0\times\cdots\times f'_{r-1} \times \gr{r,n
}\times\cdots\times\gr{n-1,n}\big)\bigcap\ff n.
$$
By the induction hypothesis the intersection above is of dimension
dim\,$\ff n-r$~ in view of the following easy

 {\bf Remark.\ }\em The validity of \rf{claim} for a given
$n$~ implies properness of the ``partial'' intersection
\vspace{-10pt}
$$
(f_0\times\cdots\times \gr{r,n+1}\times \cdots\times f_n)\cap
{\ff{n+1}},
$$
where one (or more) of the hypersurfaces $f_r\subset\gr{r,n+1}$~ is
replaced by the corresponding full grassmannian\rm.
\qed(for \rf{claim}) (Feb.2'96)

\vskip15pt
\centerline{\bf REFERENCES} \vskip2pt\parskip2pt\baselineskip13pt
\bi
\item[]{[H]  }R. Hartshorne, {\it Algebraic Geometry},
GTM \# 52 Springer--Verlag (1977).

\item[]{[K] }\ S.L. Kleiman with A. Thorup, {\it ``Intersection theory and
enumerative geometry: A decade in review'', in} Algebraic geometry: Bowdoin
1985, S. Bloch, ed., AMS Proc. of Symp. Pure Math, {\bf 46-2}, p.321-370
(1987).

\item[]{[KT]} S.~L.~Kleiman \& A.~Thorup, {\it Complete bilinear forms},
in Algebraic Geometry, Sundance, 1986, eds. A. Holme and R. Speiser,
pp. 253-320, Lect. Notes in Math. 1311, Springer--Verlag, Berlin,
(1988).

\item[]{[KTB]} $\underline{\hskip2.6cm}$, {\it A Geometric Theory of
the Buchsbaum--Rim Multiplicity},
J. Algebra \bf167\rm, 168-231 (1994).

\item[]{[L] } D. Laksov
{\it Completed quadrics and linear maps}, in Algebraic geometry:
Bowdoin 1985, S. Bloch, ed., AMS Proc. of Symp. Pure Math., {\bf
46-2}, p.371-387 (1987).

\item[]{[NT]} M. S. Narasimhan \& G. Trautmann, \it Compactification of
$M_{\p{\!3}}(0,2)$~ and Poncelet pairs of conics\rm, Pacific J. Math.
\bf145-2\rm, p.255-365
(1990).

\item[]{[P]} R. Pandharipande, {\it ``
 Notes On Kontsevich's Compactification Of The Moduli Space Of Maps},
Course notes, Univ. Chicago (1995).

\item[]{[BS]} B. Sturmfels, {\it ``Gr\"obner basis and convex
polytopes''},
Lectures notes at the Holiday Symp. at N. Mexico State Univ.,
Las Cruces (1994).

\item[]{[V]} I. Vainsencher, \it``Conics multitangent to a plane
curve''\rm, in preparation.
\ei
\vskip10pt
\parskip0pt\baselineskip10pt\obeylines
 \n Departamento de Matem\'atica
 \n Universidade Federal de Pernambuco
 \n Cidade Universit\'aria 50670--901 Recife--Pe--Brasil
 \n email: israel@dmat.ufpe.br

\enddocument